\documentclass[a4paper,leqno,11pt,twoside]{amsart}
\usepackage{amsmath,amsfonts,amssymb,amsthm,amscd}
\usepackage[utf8]{inputenc}
\usepackage[english]{babel}
\usepackage[colorlinks=true,citecolor=blue, urlcolor=blue, linkcolor=blue,pagebackref]{hyperref}
\usepackage[top=1.1in,bottom=1.1in,left=1.in,right=1.in]{geometry} 
\usepackage{graphicx}
\usepackage{wrapfig}
\usepackage{paralist}
\usepackage{tabto}
\usepackage{standalone}
\usepackage{xfrac}
\usepackage{float}
\usepackage{tikz-cd}
\usepackage{mathrsfs}

\usepackage{pgf,tikz}
\usetikzlibrary{arrows,chains,
 decorations.pathreplacing,
shapes,%
}

\usepackage[all]{xy} 

\usepackage{enumerate}
\usepackage{xcolor}
\usepackage{aliascnt}

\theoremstyle{plain}
\newtheorem{thm}{Theorem}[section]

\newtheorem{thmIntr}{Theorem}

\newaliascnt{propIntr}{thmIntr}

\newaliascnt{corIntr}{thmIntr}
\newtheorem{corIntr}[corIntr]{Corollary}

\newaliascnt{lem}{thm}
\newtheorem{lem}[lem]{Lemma}
\aliascntresetthe{lem}

\newaliascnt{claim}{thm}
\newtheorem{claim}[claim]{Claim}
\aliascntresetthe{claim}

\newaliascnt{cor}{thm}
\newtheorem{cor}[cor]{Corollary}
\aliascntresetthe{cor}

\newaliascnt{prop}{thm}

\aliascntresetthe{prop}

\theoremstyle{definition}

\newaliascnt{rem}{thm}
\newtheorem{rem}[rem]{Remark}
\aliascntresetthe{rem}

\newaliascnt{defn}{thm}
\newtheorem{defn}[defn]{Definition}
\aliascntresetthe{defn}

\newaliascnt{ex}{thm}
\newtheorem{ex}[ex]{Example}
\aliascntresetthe{ex}

\newaliascnt{setup}{thm}
\newtheorem{setup}[setup]{Set-up}
\aliascntresetthe{setup}

\numberwithin{equation}{section}

\newcommand\nc{\newcommand}
\nc\mf\mathfrak
\nc\mc\mathcal
\nc\mb\mathbb
\nc\msf\mathsf
\nc\mscr\mathscr

\def\cI{\ensuremath{\mathcal{I}}}

\def\cO{\ensuremath{\mathcal{O}}}

\DeclareMathOperator{\Hilb}{Hilb}

\definecolor{applegreen}{rgb}{0.55, 0.71, 0.0}

\title[On the locus of curves mapping to a fixed target]{On the locus of curves mapping to a fixed variety}
\author[J.~Lam, F.~Moretti and G. ~Passeri ]{Yeuk Hay Joshua Lam, Federico Moretti and Giovanni Passeri}
\address{Institut für Mathematik, Humboldt-Universität zu Berlin, Unter den Linden 6, 10099 Berlin, Germany}
\email{joshua.lam@hu-berlin.de} 

\address{Department of Mathematics, State University of New York at Stony Brook, 100 Nicolls Rd, Stony Brook, NY 11794, USA}
\email{federico.moretti@stonybrook.edu}
\email{giovanni.passeri@stonybrook.edu} 

\begin{document}

\begin{abstract}
Suppose $Y$ is a smooth variety equipped with a top form. We prove a simple theorem giving a sharp lower bound on the geometric genus of a family of subvarieties of  $Y$,  in terms of the dimension of this family. 
As an application,  we give  upper bounds on the dimension of the locus of curves admitting a map to $Y$, when  $Y$ is a hypersurface or an abelian variety. In the latter case we characterize the  families of curves attaining this bound as the families of degree $2$ branched covers of a  fixed curve.

\end{abstract}


\maketitle

\setcounter{tocdepth}{1}
\section{Introduction}

\quad The goal of the present paper is to study the geometric genus of families of subvarieties on any smooth variety $Y$ with positive canonical bundle. In particular given a smooth variety $Y$ we will give upper bounds on the dimension of the locus $M_{g,Y}^{bir}\subset \mathcal M_g(\mathbb C)$ defined as follows \begin{eqnarray*}
    M_{g,Y}^{bir}:={\{[C] \in \mathcal M_g \ | \ \exists f:C \to Y \text{ such that $f$ is birational onto its image}\}} \subset \mathcal M_g
    \end{eqnarray*}

    These loci are countable unions of closed subvarieties (c.f. \autoref{set up hypersurfaces}) so it makes sense to talk about their dimension.
We first give the following simple-to-state application of our main result.

\begin{thmIntr}
   \label{hypint}
    Let $Y\subset \mathbb P^{n+1}_\mathbb{C}$ be a very general hypersurface of degree $d \ge 2n-1$ (if $n\ge 4$) or of degree $d \ge 2n+1$ (if $n=2,3$), then  \[\dim M_{g,Y}^{bir} \le g-(d-2n).\]

  In particular, for a very general curve $C$ of genus $g$: 
    \begin{enumerate}
        \item if $d\ge 2n$, any map $C\to Y$ is constant;
        \item if $d=2n-1$ and $g\ge 3$, any map $C\to Y$ factors through a rational curve. 
    \end{enumerate}
 
\end{thmIntr}

\quad  Let us remark that \autoref{hypint} is almost optimal:  any hypersurface $Y$ of degree $<2n$ contains a line, hence for any curve $C$ we have non-constant maps $C\to Y$. For $n=2$ any surface of degree $2n=4$ contains a rational curve. Whether the bound on the dimension is sharp is a difficult question (and the answer is most likely no).

\quad \autoref{hypint} above is obtained by combining pioneering results of Ein \cite{Einhypersf} and Voisin \cite{Voisinhypersf} on the positivity of the canonical bundle of subvarieties of general hypersurfaces, and   the following general theorem. 
\newpage
\begin{thmIntr}\label{estimate}
    Let $Y$ be a smooth variety over an algebraically closed field with $h^0(\omega_Y)\ge 1$ and $B\subset \mathrm{Hilb}_Y$ be an irreducible reduced subscheme parametrizing (possibly singular) irreducible reduced subvarieties of codimension $r$ and geometric genus $g$. Let\footnote{We mean that for a general closed point $b\in B$, the fiber $H_b$ is smooth.  Such a desingularization always exists over a field of charachteristic zero.}  $f:\mathcal{H}\to B$ be a relative desingularization of the universal family, and write $H_b$ for the fiber at $b\in B$.  Suppose that the morphism $\pi:\mathcal H \to Y$ is generically smooth. Then,   for any $b\in B$,
   \[ g\ge h^0({\omega_Y}_{|_{\pi({H_b})}})+\mathrm{dim}(B)-r.\]

\end{thmIntr}
\indent The bound is sharp and all the hypotheses are necessary, see  \autoref{section:examples} for a discussion. There is an analogous statement for infinitesimal families, see \autoref{infinitesimal}. One of the key points of \autoref{estimate} is that we do not put any restrictions on the type of singularities that the general fiber $H_b$ has. 

In the special case of non degenerate curves in an abelian variety $A$, we obtain stronger results; this improvement comes from an analysis of  the associated \emph{Gauss maps}. The case of curves with a map to $A$ has attracted considerable recent attention, for example through the Putman--Wieland conjecture and its recasting by Landesman--Litt in algebro-geometric terms: see  \autoref{section: motivation} for further discussions. Recall that given a family of curves $f:\mathcal C \to B$ the isotrivial factor of the family of Jacobians is the largest abelian subvariety $A_f$ containted in the Jacobian of all the fibers $\mathrm{Jac}(C_b)$. Up to compactifying $\mathcal C$ and base change, $A_f$ is the quotient of  $\mathrm{Pic}^0(B)\to \mathrm{Pic}^0(\mathcal C)$.
For example, we obtain (see also \autoref{curvesinab} for the case of curves mapping birationally to $A$)

\begin{thmIntr}\label{thm:genus-bounds}
Let $f:\mathcal{C}\to B$ be a  family of curves of genus $g$ over $\mathbb C$, whose associated map $B\rightarrow \mc{M}_g$ is finite onto its image. Let $A_f$ of dimension $q_f\ge 2$ be the isotrivial isogeny factor of the associated family of Jacobians. Then,
\[g\geq 
\begin{cases}
 \mathrm{dim}(B)+q_f  &\textrm{if}    \quad\mathrm{dim}(B)\le \frac{1}{2}(q_f-1)\quad \mathrm{or} \quad  2q_f-5\le \dim B \le 2q_f-2 \\
 \frac{3 \mathrm{max}\{\dim B+1, q_f-1\}}{2} +1   &\textrm{if} \quad \, \,\frac{1}{2}(q_f-1) \le  \mathrm{dim}(B) \le  2q_f-5 \\
 2q_f-1+\frac{\dim B}{2}&\textrm{if}  \quad \, \, 2q_f-2\leq \mathrm{dim}(B).
\end{cases}
\]
Moreover, suppose 
$\mathrm{dim}(B)>2q_f-2$, then the above bound is sharp and the only family realizing this bound is the family of degree $2$ covers of a curve of genus $q_f$ branched over $n=\mathrm{dim}(B)$ points (hence $B$ is a generically finite cover of $C^{(n)}$ the $n$-fold symmetric product of $C$ and $A_f$ is isogeneous to $\mathrm{Jac}(C)$).
\label{iso}
\end{thmIntr}
The second part of the theorem can be rephrased in various ways. For example, the following has  a more topological flavor:
\begin{corIntr}
    Let $\mathcal C \to B$ be a genus $g$ fibration such that the map  $B\rightarrow \mathcal M_g$ is finite onto its image and $\mathrm{dim}(B)\ge \frac{2g-2}{3}$. Then
    \[
    \mathrm{rk}(H^0(\Omega^1_\mathcal C)\to H^0(\omega_{C_b}))\le \bigg\lfloor \frac{g}{2}+\frac{1}{2}-\frac{\mathrm{dim}(B)}{4} \bigg\rfloor,
    \]
    for any $b\in B.$
    The bound is sharp, moreover, if $\mathrm{dim}(B)>\frac{2g-2}{3} $ equality occurs if and only if we are in the following situation: there exists  an integer $k\geq 1$ and a genus $\big\lfloor\frac{g}{2}+\frac{1}{2}-\frac{\mathrm{dim}(B)}{4}\big\rfloor$ curve $C$ such that $\mathcal{C}\rightarrow B$ is a family of degree 2 covers branched at $k$ points of $C$.
    \label{corrrr}
\end{corIntr}

 Given an abelian variety $A$, let us denote by
\[
 M_{g,A}^{nd}:={\{[C] \in \mathcal M_g \ | \ \exists f:C \to A \text{ such that $f$ is non degenerate}\}} \subset \mathcal M_g.
\]
With this notation we obtain
\begin{corIntr}\label{cor: abelian}
 Let $A$ be an abelian variety of dimension $q$ and consider $g> 2q-2$, then
\[
 {\mathrm{dim}(M_{g,A}^{nd})}\le 2g-4q+2. 
\]

The bound is sharp, and equality occurs if and only if we are in the following situation: there exists  an integer $k\geq 1$ and a genus $q$ curve $C$ such that $A$ is isogenous to  $\mathrm{Jac}(C)$, and $\mathcal{C}\rightarrow B$ is the family of degree 2 covers branched at $k$ points of $C$.
\label{coriso}
\end{corIntr}

Let us stress that in \autoref{iso}, \autoref{corrrr} and \autoref{coriso} we are not assuming that the map $C_b\to A$ is birational onto its image (as opposed to \autoref{hypint}, where we do assume this). The motivation here is  of a somewhat different and more Hodge-theoretic nature: we are interested in bounding the dimension of the fixed part of  families of Jacobians. 

 In order to illustrate the bounds that we get, we collect the known results in the case $q_f=3$ (and $q=3$ in the notation of  \autoref{cor: abelian}), i.e. the case of abelian threefolds, into the following table. As far as we are aware, these are the optimal known results for the question we are interested in (and the same is true for higher dimensional $A$). 
\begin{center}

\begin{tabular}{ |p{2cm}||p{7,7cm}|p{4,7cm}|  }
				\hline
				\multicolumn{3}{|c|}{Families of Jacobians with an abelian threefold $A$ as an isotrivial factor} \\
				\hline
			genus $g$ & Dimension of the biggest component of $ M_{g,A}^{nd}$ & Sharpness\\ 
				\hline
			
			 3 & \quad 0 (Torelli Theorem)  & for any $A$ \ \\
			4 &  \quad 1  (Torelli Theorem)  & for some $A$ \cite{MR1179333}   \\
		 5 &  \quad $\le 2$ (\autoref{iso}) & ??  \\
				6 & \quad  $\le 3$ (\autoref{iso}) & ?? \\
				 $g\ge 7$    & \quad $2g-10$ (\autoref{iso}) & for any $A$  \\
				\hline
			\end{tabular}

\end{center}

\subsection{Motivations}\label{section: motivation} We now give some  motivation behind the results stated  above. The study of curves and subvarieties of hypersurfaces received considerable attention since the work of Clemens on curves in hypersurfaces \cite{clemens} (for a comparison of his Theorem with our result see \autoref{cle:curves vs hyp} and the following comments), followed by the aforementioned works by Ein and Voisin \cite{Einhypersf, Voisinhypersf}. More refined statements were obtained in case of divisors by Xu \cite{xu} (for a different perspective and variations of this problem see also \cite{lopez, clemensran}). In recent years the positivity results  obtained found remarkable applications in the study of rational maps $Y\dashrightarrow \mathbb P^n$. Recall that the $\textrm{degree of irrationality}$ of a variety $Y$ is the minimum degree of a dominant rational map to a projective space of the same dimension. Bastianelli, De Poi, Ein, Lazarsfeld and Ullery proved in \cite{BDELU} that the degree of irrationality $\mathrm{irr}(Y)$ of a very general hypersurface $Y$ of degree $d \ge 2n+1$ is $d-1$. In a somehow similar fashion Lazarsfeld and Martin in \cite{roboli} extended this study to correspondences between hypersurfaces of the same dimension, proving that for $Y_1,Y_2\subset \mathbb P^{n+1}$ very general hypersurfaces of degrees $d_1,d_2\ge 2n+1$, then any non degenerate correspondence $Z\subset Y_1\times Y_2$ satisfies $\mathrm{deg}(Z\to Y_i)\ge d_i-1$. In \cite{Passeri} the third author extends this study to stable correspondences between hypersurfaces of possibly different dimensions, \autoref{hypint} (and more specifically \autoref{hypersf} in Section 5) plays an important role since if there is a non constant map $X\to Y$ for very general hypersurfaces one can construct a special stable correspondence between them (for more details see \cite[Example 2.2]{Passeri}).

\quad The study of families of Jacobians with an isotrivial factor has a long history. The literature focuses  on extreme situations: either families of dimension $1$ or families over the whole $\mathcal M_{g,n}$. The first case falls in the field of research of irregular surfaces: the first inequality in this direction was proven by Xiao in  \cite{xiao}, proving that if $B=\mathbb P^1$ then $g\ge 2q_f-1$ (in the notation of \autoref{iso}), this fails if the base is not $\mathbb P^1$ by  constructions of Pirola and Albano--Pirola \cite{MR1179333, AlbPir}; under some generality  assumptions on the Clifford index  of the general curve of the family a slightly weaker inequality holds by work of Barja, González Alonso and Naranjo \cite{xiaobgn}. General results in this direction seem rather difficult to come by. As far as we know, the strongest result is due to Xiao \cite[Corollary 3]{gang1987fibered}, who proves that $g\geq \frac{6(q_f-1)}{5}$ for families of curves over an arbitrary one dimensional base. It is perhaps worth mentioning that this gives a partial answer to a question of Landesman--Litt \cite[Question 5.16]{LL2}, which seems to have been overlooked. 

\quad At the opposite extreme, one can study  the case when  the base is a finite cover of $\mathcal{M}_{g,n}$. One notable case is when   the family of curves is a $G-$cover of a general curve of type $(g, n)$, for some finite group $G$ prescribing the monodromy. The Putman--Wieland conjecture \cite{PW} predicts that in this situation the family of Jacobians has  no  isotrivial factor. The conjecture is known to be false for $g=2$ by work of Markovìc \cite{MR4549123}, for a short proof see also \cite{LLPrill}. Work of Landesman and Litt settles the conjecture for $G$ of cardinality $\le g^2$, \cite{litt2022canonical}. We refer the reader to \cite[Section 5.10]{LL2} for a   comprehensive introduction, as well as  some related interesting open questions.

\textbf{Proofs.} Let us give a brief indication of the proof of \autoref{estimate},  carried out in \autoref{section:bound}. The key is to study, in the setup of \autoref{estimate},  the contraction map  \[
   c: \bigwedge^r T_bB \otimes H^0({\omega_Y}_{|_{\pi(H_b)}})\to H^0(\omega_{H_b});
    \]
    roughly, this map takes a top form of $Y$, and contracts it with $r$ independent directions in the base $B$, thus producing a top form on fibers. In the case $\pi(H_b)$ is smooth this contraction map is induced by the adjunction formula. The bound in \autoref{estimate} is then obtained by giving a lower bound of the rank of this contraction map. This is achieved by finding a suitable (non-canonical) \emph{volume decomposition} $T_bB=\bigoplus V_i$ such that the contraction is never zero on tensors of the type $\omega\otimes (\partial_{t_1}\wedge \dots \wedge\partial_{t_{n-1}})$, with $\partial_{t_i}\in V_i$ non-zero. The bound then follows by Hopf's lemma (c.f. \autoref{hopf}). This volume decomposition is obtained by pulling back a filtration of $\pi^*\mathcal{N}_{\pi(H_b)|Y}$ via $d\pi:\mathcal N_{H_b|Y}=T_bB\otimes \mathcal O_{H_b}\to \pi^*\mathcal{N}_{\pi(H_b)|Y}$ (for more details see \autoref{rank}). When the normal bundle is the dual of the \emph{kernel}\footnote{We say that a vector bundle is \emph{kernel sheaf} of a morphism (or a rational map) $X\to \mathbb PV$ if it is the kernel of the evaluation $V^\vee\otimes \mathcal O_X\to \mathcal O_X(1)$.} of some non degenerate rational map to $\mathbb P^r$ it is possible to obtain better bounds (see section \autoref{23}).  It is perhaps worth pointing  out that our  proof is elementary and completely self contained, the results in this section hold over any algebraically closed field.

    \quad In \autoref{section: curves isoztivial}, \autoref{iso} is obtained using the improved version of \autoref{estimate} together with some bounds of the degree of $\emph{kernel}$ bundles, see \autoref{lemma:degree3} and \autoref{degreeez}. In this case the normal bundle of the curves in abelian variety is the dual of the \emph{kernel} bundle of the Gauss map. The results in this section hold over $\mathbb C$.

    \quad In the last section \autoref{hypint} is obtained applying \autoref{estimate} together with the aforementioned results of Ein and Voisin that provides a lower bound for $h^0({\omega_Y}_{|_{\pi({H_b})}})$. This section is independent of \autoref{section: curves isoztivial}. The results in this section hold over $\mathbb C.$

\textbf{Acknowledgements.} We would like to thank Gavril Farkas and Robert Lazarsfeld for mentoring and interesting discussions. We thank Victor González Alonso and Pietro Pirola for useful conversations about contractions. We are grateful to Raju Krishnamoorthy for organizing a reading group on the work of Landesman and Litt. We would like also to acknowledge Rita Pardini and Pietro Pirola for some discussions on Xiao's conjecture and, in particular, for showing us the paper \cite{mpp}. Finally we thank Andrea Di Lorenzo for some discussions on the Hilbert scheme.

\quad The first author is supported by a Dirichlet Fellowhsip at Humboldt Universit\"at zu Berlin. The second author is supported by the ERC Advanced Grant SYZYGY. This project has received funding from the European Research Council (ERC) under the European Union Horizon 2020 research and innovation program (grant agreement No. 834172). Part of this project has been carried out during the visit of the second author at the State University of New York at Stony Brook funded by the Berlin Mathematical school.

\section{Geometric genus of subvarieties covering a smooth variety with a top form}
\label{section:bound}
\quad The goal of this section is to prove \autoref{estimate}. We begin with some preliminaries.
\subsection{Preliminaries}
We will use the following lemma in linear algebra, usually attributed to Hopf.
\begin{lem}[Hopf's lemma]
Let $V_1,\dots,V_n,W$ be vector spaces of dimension $l_1,\dots,l_n,m$  over an algebraically closed field and $\Psi :\bigotimes_{i=1}^nV_i \to W$ be a linear map such that $\Psi(v_1\otimes \dots \otimes v_n)\neq 0$ for any $v_1,\dots, v_n\neq 0$ then $m\ge \sum_i l_i-n+1$.
\label{hopf}
\end{lem}

\quad The following is a slightly stronger version of \cite[Proposition 4.5]{mpp}, we also present a variation of the proof in \cite{mpp}.
\begin{lem}
 Let $X$ be an irreducible generically reduced variety over a field $ \mathbb K=\overline{\mathbb K}$, $\mathcal P\in \mathrm{Coh}(X)$  a coherent sheaf of generic rank $r$ and $V\otimes \cO_X\to \mathcal P$ be a morphism of generic rank $r$ where $V$ is finite dimensional vector space. Then there is a decomposition $\bigoplus_{i=0}^r V_i=V$ with $V_i\neq 0$ for $i\ge 1$ such that for any  non zero $v_i\in V_i$
the map\[
\langle v_1,\dots,v_r \rangle \otimes \cO_X \to \mathcal P
\]
is of generic rank $r$. In particular 
\[
\mathrm{rank}_{\mathbb K}(\mathrm{det}_U)\ge\mathrm{dim}(V)-\mathrm{dim}(V_0)-r+1, \]
where $\mathrm{det}_U: \bigwedge^r V\to H^0((\bigwedge^r \mathcal P)_{|_U}) $ for any open $U\subset X.$
\label{rank}
\end{lem}
\begin{proof}
First we realize a (non canonical) volume filtration $W_0\subsetneq W_1\subsetneq\dots \subsetneq W_r=V$. We set $W_0$ as a maximal subspace such that the map $W_0\otimes \cO_X\to \mathcal P$ is of generic rank $0$ (notice that if $\mathcal P$ has torsion this does not necessarily mean that the map is zero). Then inductively we set $W_{i-1}\subsetneq W_i$ as a maximal subspace such that $W_i\otimes \cO_X\to \mathcal P$ is of generic rank $i$. Notice that we have $W_i=V$ if and only if $i=r$ since the morphism $V\otimes \cO_X\to \mathcal P$ is of generic rank $r$. The (non canonical) decomposition will be $V_0=W_0$ and $V_i$ a complement of $W_{i-1}$ in $W_i$ (i.e. $V_i$ is some subspace satisfying $W_i=W_{i-1}\oplus V_i$).

\quad For the second part notice that $\mathrm{det}_U(v_1\wedge \dots \wedge v_r)=0$ if and only if the morphism \[
\langle v_1,\dots,v_r\rangle \otimes \cO_X \to\mathcal P
\]
is of generic rank $<r$. Now if we pick non zero $v_i\in V_i$ for $i=1,\dots, r$ and the above map is of generic rank $<r$ then there is a minimal $k$ such that 
\[
\langle v_1,\dots,v_k\rangle \otimes \cO_X \to \mathcal P
\]
is of generic rank $k-1$, but then $(\langle v_k \rangle \oplus W_{k-1})\otimes \cO_X \to\mathcal P$ is of generic rank $k-1$. We get $v_k\in W_{k}$, contradiction. The bound on the rank follows directly by Hopf's lemma applied to the restriction of the determinant $\bigotimes V_i \to  H^0((\bigwedge^r \mathcal P)_{|_U})$.

\end{proof}
\begin{rem}
\quad Let us denote with $\mathcal P_k$ to be the torsion free part of the quotient of $\mathcal P$ by the image of $W_{k-1}\otimes \cO_X\to \mathcal P$ and by $Z(V_i)$ the zero locus of the map $V_i\otimes \cO_X\to \mathcal P_i$. For any  open $U\subset X$ such that $\mathcal 
P_{|_U}$ is locally free we get that the image of $\mathrm{det}_U: \bigotimes_{i=1}^r V_i\to H^0(\bigwedge^r \mathcal P)$ factors through $H^0(\bigwedge ^r \mathcal P_{|_U}\otimes \cI_{\cup_{i=1}^r Z(V_i)})$. This will be particularly useful in \autoref{23}. \label{remardeterminant}
\end{rem}
\begin{setup}
Let $Y$ be a smooth variety of dimension $n$ over an algebraically closed field with one canonical form. Let $B\subset \mathrm{Hilb}_Y$ be an irreducible reduced subscheme parametrizing irreducible reduced subvarieties. Let us consider a relative smoothing $f:\mathcal{H}\to B$ of the universal family (i.e. for the general closed point $b\in B$ we have $H_b$ smooth). Suppose that the general $H_b$ is  of geometric genus $p_g(H_b)\le g$ and that the map $\pi:\mathcal{H}\to Y$ is generically smooth with $\pi(H_b)$ of codimension $r$. Moreover let $U_b\subset H_b$ be an open subvariety such that $\pi:U_b\to\pi(U_b)$ is an isomorphism.
\end{setup}
 \quad Let us fix a general $b\in B$. We have a natural contraction map \[
   c: \bigwedge^r T_bB \otimes H^0({\omega_Y}_{|_{\pi(H_b)}})\to H^0(\omega_{H_b}),
    \]
    where $c(\partial_{t_1}\wedge\dots\wedge\partial_{t_r}\otimes \omega)_p=\pi^*\omega(\partial_{t_1},\dots,\partial_{t_r},-)_{|_{(T_{p}H_b)^{n-r}}}$ where we abuse notation and denote with $\partial_{t_i}$ a local lifting of $\partial_{t_i}$ to a small neighborhood of $p\in H_b$ (a local lifting is well defined up to $T_pH_b$ hence the restriction of the form is independent of the chosen local lifting). We can now prove.
\begin{thm}

In the notation as above there is a decomposition $T_bB=\bigoplus_{i=1}^r V_i$ such that for any non zero $v_i\in V_i,\omega\in H^0({\omega_Y}_{|_{H_b}})$ the contraction
\[
c( v_1\wedge\dots\wedge v_r\otimes \omega)\in H^0(\omega_{H_b})
\]
is non zero. In particular  
\[g\ge h^0({\omega_Y}_{|_{H_b}})+\mathrm{dim}(B)-r.\]

\label{main}

\end{thm}

\begin{proof}

If $\omega\in H^0(\omega_Y)$ is generically non zero on $\pi(H_b)$ then $c(\partial_{t_1}\wedge\dots\wedge\partial_{t_r}\otimes \omega)=0$ if and only if for the general $p\in H_b$ we have 
    \[
    \langle d\pi_p(\partial_{t_1}),\dots,d\pi_p(\partial_{t_r}),T_{\pi (p)}\pi(H_b)\rangle \subsetneq T_{\pi(p)}Y.
    \]
   \quad  Let us write the condition in terms of normal sheaves, on $ H_b$ we have a commutative diagram with exact rows\[
    \begin{tikzcd}
    0 \arrow{r} & \mathcal T_{H_b}\arrow{r}& {\mathcal T_{\mathcal H}}_{|_{H_b}}\arrow{r} \arrow{d}{d\pi}& \mathcal N_{H_b/\mathcal H}=T_bB\otimes \mathcal O_{H_b}\arrow{d}{d\pi}\arrow{r} & 0 \\
    &  & {\pi^*\mathcal T_{Y}}_{|_{\pi(H_b)}}\arrow{r}& \mathcal P \arrow{r} & 0,
    \end{tikzcd}
    \]
    where  $\mathcal P $ is the torsion free part of the cokernel of $d\pi:\mathcal T_{H_b}\to{\pi^*\mathcal T_{Y}}_{|_{\pi(H_b)}}$. 
     The vertical lines are induced by the differential. Notice that on $U_b$ the above diagram becomes
    
    \[    \begin{tikzcd}
    0 \arrow{r} & \mathcal T_{U_b}\arrow{r} \arrow{d}{d\pi}& {\mathcal T_{\mathcal H}}_{|_{U_b}}\arrow{r} \arrow{d}{d\pi}& \mathcal N_{U_b/\mathcal H}=T_bB\otimes \mathcal O_{U_b}\arrow{d}{d\pi}\arrow{r} & 0 \\
     0 \arrow{r} & \mathcal \pi^*\mathcal T_{\pi(U_b)}\arrow{r} & {\pi^*\mathcal T_{Y}}_{|_{U_b}}\arrow{r}& \mathcal \pi^*\mathcal N_{\pi(U_b)/Y}\arrow{r} & 0.
    \end{tikzcd}
    \]
    The above condition $c(\partial_{t_1}\wedge\dots\wedge\partial_{t_r}\otimes \omega)=0$ reads that \[\begin{tikzcd}
    \langle \partial_{t_1},\dots,\partial_{t_r}\rangle \otimes \mathcal O_{H_b} \arrow{r}{d\pi} & \mathcal P
    \end{tikzcd}
    \] 
    is not of generic rank $r$. Now since the map is generically smooth we get that ${\mathcal T_{\mathcal H}}_{|_{H_b}}\to {\pi^*\mathcal T_{Y}}_{|_{\pi(H_b)}}$ is generically surjective. We deduce that $T_bB\otimes \mathcal \cO_{H_b}\to \mathcal P$ is of generic rank $r$ (i.e. generically surjective).   Now by \autoref{rank} we can find a decomposition $T_bB=\bigoplus_{i=0}^r V_i$ such that  for any  non zero $v_i\in V_i$
the map\[
\langle v_1,\dots,v_r \rangle \otimes \cO_X \to \mathcal P
\]
is of generic rank $r$. The induced map $\bigotimes_{i=1}^r V_i\otimes H^0({\omega_Y}_{|_{H_b}}) \to H^0(\omega_{(H_b)})$ is non zero on simple tensors.
If $H_b$ is smooth then any $v_0\in T_bB$ has volume, meaning that  $V_0=0$ or equivalently that $d\pi (v_0) \in H^0(\mathcal P)$ is 
 generically non zero for any $v_0 \in T_bB$ as the following argument shows. Consider the infinitesimal deformation in the direction of $v_0$
 \[
 \begin{tikzcd}
 H_{v_0} \arrow{d}\arrow{r} & Y \\
 \mathrm{Spec} \, \mathbb K[\epsilon] \arrow{r}{v_0} & B\subset \mathrm{Hilb}_Y.
 \end{tikzcd}
 \] The condition $d\pi(v_0)=0$ implies $\pi(U_{v_0})=\pi(U_b)$ scheme theoretically. We deduce that the image of $v_0$ in $T_{\pi(H_b)}(\mathrm{Hilb}_Y)=\mathrm{Hom}_{\pi(H_b)}(\mathcal{I}_{\pi(H_b)}/\mathcal I_{\pi(H_b)}^2,\mathcal O_{\pi(H_b)})$ is generically zero, which implies that it is zero (it is a generically zero section of a torsion free sheaf). This implies that $v_0$ lies in the kernel of $T_bB\to T_{\pi(H_b)}\mathrm{Hilb}_Y$, that contradicts the inclusion $B\subset \mathrm{Hilb}_Y$.

 Recall that the induced map $\bigotimes_{i=1}^r V_i\otimes H^0({\omega_Y}_{|_{H_b}}) \to H^0(\omega_{(H_b)})$ is non zero on simple tensors. The bound on the rank follows from Hopf's lemma together with $V_0=0$.

\end{proof}
\quad As mentioned in the introduction the above proof can be adapted to prove the following version that works also  for infinitesimal families.

\begin{cor}
    Let $Y$ be a scheme and $B\subset \mathrm{Hilb}_Y$ be an infinitesimal subscheme with closed point $b$ \footnote{equivalently a fat point}, where $\Hilb_Y$ is the Hilbert scheme 
 parametrizing (possibly singular) irreducible reduced subschemes of $Y$ of dimension $n-r$ and geometric genus $g$.  Let  $f:\mathcal{H}\to B$ be a relative desingularization of the universal family, suppose that the morphism $\pi:\mathcal H \to Y$ has generically surjective differential. Suppose moreover that ${\Omega_Y}_{|_\pi(H_{b})}$ is locally free of rank $n$ and that $h^0({ \bigwedge ^n \Omega_Y}_{|_{\pi(H_b)}})\ge 1$. Then the following holds 
   \[ g\ge h^0\Big({\bigwedge^n \Omega_Y}_{|_{\pi(H_b)}}\Big)+\mathrm{dim}(T_bB)-r.\]
   \label{infinitesimal}
\end{cor}
Let us remark that by relative desingularization we mean that $H_b$ is smooth for the closed point $b\in B$. This means that we can consider only the deformations  $T_bB\subset T_{\pi(H_b)}(\mathrm{Hilb}_Y)$ that are preserving the geometric genus.

In particular if $\pi({H_{b}})$ is smooth we can deduce an upper bound on the sections of the normal bundle.
\begin{cor}
Under the same hypotheses as \autoref{infinitesimal}, suppose moreover that $\pi(H_{b})$ is smooth. Then \[
h^0(\mathcal N _{\pi(H_{b})|Y})\le g-h^0\Big({\omega_Y}_{|_{\pi(H_b)}}\Big)+r.
\]
\end{cor}
\subsection{Examples and discussions}\label{section:examples}
We discuss now the hypotheses in \autoref{estimate}, prove that the bound is sharp, as well as to  provide some elementary examples.

\textbf{Generic smoothness and relative desingularization.}
Notice that generic smoothness in characteristic zero reduces to $\mathcal H \to Y$ dominant. This is necessary, for instance there are hypersurfaces $Y$ of high degrees containing planes, any curve $C$ has a non constant map $C\to Y$ factoring through the plane. In the hypothesis of \autoref{estimate} it is in principle possible that generic smoothness is implied by being dominant even in positive characteristic (this is true, for instance, when $\pi(H_b)$ is of codimension $1$). What can fail in positive characteristic is the existence of a relative desingularization. This is the case of quasi-elliptic surfaces, that is smooth surfaces $S\to B$ in characteristic $2,3$ where $B$ is smooth and all the fibers are cuspidal rational curves, in characteristic $3$ the local equation is $t=y^2-x^3$ (c.f. \cite{bombierimumford}). In this case, since the total family is already smooth over the base field, a relative desingularization over $B$ cannot exist. Another example is given by Fermat hypersurfaces, see \cite{shiodafermat}.

\textbf{Irreducibility.} The irreducibility condition is also necessary as the following example shows.
\begin{ex}
Let $A=E\times S$ be an abelian threefold, with $S$  abelian surface and $E$ an elliptic curve. Take $C\subset S\times {p}$ with $p\in E$ a curve of genus $2$. Then the family $\mathcal C\to A\times A$ where $C_{a,b}=\{C+a\}\cup\{C+b\}$ is a $6$ dimensional family of curves of genus $4$ contradicting the bound in the above theorem.
\end{ex}
If $r=1$ the irreducibility condition can be avoided.

\begin{ex}[Sharpness of the bound] Severi varieties provides example of family of genus $g$ curves moving in a family of dimension $g$ on a $K3$ surface or on an abelian surface \cite{chen}\cite{severi}. In case $r=1$, they provide an example where the bound of \autoref{estimate} is sharp. In case $r\ge2$ consider $S$ an abelian surface with a $g-$dimensional family of curves of genus $g$ over $B$, take any abelian variety $A$. Translating we get a $\mathrm{dim}(A)+g$ dimensional family of curves of genus $g$ over $A\times B$ in $A\times S$; this proves that the bound is sharp even for higher codimension. Notice that all this curves are \emph{degenerate} in $A\times S$ (so there is no contradiction with \autoref{iso}). A similar example can be produced in the Hilbert scheme of $n$ points $S^{[n]}$ of a $K3$ surface $S$.
\end{ex}

\subsection{Improved bounds in geometric situations}
\label{23}
Before moving to the applications mentioned in the introduction we would like to give an improved version of \autoref{rank} when the sheaf $\mathcal P^\vee$ is the \emph{kernel} sheaf of some rational map. We will then observe that in some concrete situations the normal bundle fits in this hypothesis.
\begin{lem}\label{lemma:gaussmap}
  Let $X$ be a smooth projective variety of  and consider $f:X \dashrightarrow \mathbb P^{r+1}$ a non degenerate map  such that for the general point $p\in f(X)$ we have that either $f^{-1}(p)$ is of codimension $\ge 2$ or $h^0(\mathcal O(f^{-1}(P)))=1$.
  Let $\mathcal{P}^\vee\in \mathrm{Coh}(X)$ denote the kernel sheaf of $f$\footnote{It is just a reflexive sheaf in general}, defined as 
  \[
  0 \to \mathcal P^\vee \to H^0(\mathcal O_{\mathbb P^{r+1}}(1)) \otimes \mathcal O_X \to \mathcal O_X(1)\otimes \mc{I}\to 0,
  \]
  then \[
  \mathrm{rank}(\bigwedge ^r H^0(\mathcal P) \to H^0(\mathcal O_X(1))) \ge h^0(\mathcal P).
  \]
\end{lem}
Before moving onto the proof let us mention the following  elementary fact that we will be using sometimes implicitly.
\begin{rem}
    In the notation as in the above lemma the zero locus of a section $[s] \in \mathbb P H^0(\mathcal O_{\mathbb P^{r+1}}(1) )^{\vee}\subset \mathbb PH^0(\mathcal P)$ contains $f^{-1}([s])$ and is is contained in $f^{-1}([s])\cup \mathrm{Bl}(f)$. Here $\mathrm{Bl}(f)$ denote the base locus of the map $f$ (which is in codimension $\ge 2$). 
\end{rem}
\begin{proof}We argue by induction on $r$. Notice that in case that $r=0$ the hypothesis are not satisfied but the conclusion of the lemma holds, so we may use this as our base case of the induction.
     Pick $p_1\in X$ general, and let $v_1\in H^0(\mc{O}_{\mb{P}^{r+1}}(1))^{\vee}\subset H^0(\mc{P})$ be the section corresponding   to  $f(p_1).$ The induced map $\mathcal O_X\to \mathcal P$ is vanishing along $Z(v_1)$. From  our hypothesis we deduce that 
     \begin{itemize}
         \item either $Z(v_1)=Z(f^{-1}f(p_1))$ is a non movable divisor if $X$ is a curve,
         \item or $Z(v_1)$ is of codimension $\ge 2$ if  $X$ is of dimension $\ge 2$.
     \end{itemize} 
     This implies that there cannot be any other section $s\in H^0(\mathcal P)$ such that $s\wedge v_1=0$, otherwise $Z(v_1)$ would be a moving divisor \footnote{If one wants to apply the strategy of \autoref{rank} this means that $V_1=\langle v_1 \rangle$} (the $2$ sections $s,v$ would span a subline bundle of $\mathcal P$, hence the zero locus of $v_1$ would contain a moving divisor). This means that $H^0(\mathcal P)/\langle v_1 \rangle \subset H^0((\mathcal P/\mc O_X \cdot v_1)^{\vee \vee})$. But, now $\mathcal P_1=(\mathcal P/\mc O_X \cdot v_1)^{\vee \vee}$ is the dual of the kernel of $H^0(\mathcal O_{\mathbb P^{r+1}}(1))(-p_1)\otimes \mathcal{O}_X\to \mathcal O_X(1)\otimes \mathcal I_{Z(v_1)}$. Indeed we have a commutative diagram 
     \[
     \begin{tikzcd}
      0 \arrow{r} & \mathcal P_1^\vee  \arrow{r} \arrow[d,hook] &  H^0(\mathcal O_{\mathbb P^{r+1}}(1))(-p_1)\otimes \mathcal{O}_X \arrow[d,hook] \arrow{r} & \mathcal O_X(1)\arrow{d} \\
       0 \arrow{r} & \mathcal P_1^\vee  \arrow{r} \arrow{d}{v_1^\vee} &  H^0(\mathcal O_{\mathbb P^{r+1}}(1))(-p_1)\otimes \mathcal{O}_X \arrow{d} \arrow{r} & \mathcal O_X(1)\\  & \mathcal O_X \arrow{r} & \mathcal O_X &
 \end{tikzcd}
 \]
     Moreover if the point $p_1$ is generic we may suppose that either the hypothesis of the lemma are still verified or that  $h^0(\mathcal O_{\mathbb P^{r+1}}(1))(-p_1))=2$ (the base case of the induction). Indeed if $X$ is a curve the projection of $f(X)$ from $f(p_1)$ is birational onto its image as long as $h^0(\mathcal O_{\mathbb P^{r+1}}(1))(-p_1))>2$ and $p_1$ is general (hence the preimages of general points under the map $X \to \mathbb P  (H^0(\mathcal O_{\mathbb P^{r+1}}(1))(-p_1))$ are still non movable divisors). The case of $X$ of dimension $\ge 2$ is similar. 
     
     Hence we get by induction
 \[
  \mathrm{rank}(\bigwedge ^{r-1} H^0((\mathcal P_1) \to H^0(\mathcal O_X(1)))\ge h^0((\mathcal P_1)\ge h^0(\mathcal P)-1.
\]
Now the image of the original determinant contains the image of the multiplication map
\[
\begin{tikzcd}
\langle v_1 \rangle \otimes \bigwedge ^{r-1} H^0((\mathcal P_1) \arrow{r}{m_1}\arrow[d,hook] &H^0(\mathcal O_X(1)) \\
\bigwedge ^r H^0(\mathcal P) \arrow{ur}{\mathrm{det}} &
\end{tikzcd}
\]
Now since the section $v_1$ is vanishing at the  general point $p_1\in X$ we get that $m_1$ factors through $H^0(\mathcal O_X(1)\otimes \mathcal I_{p_1})$ which is strictly contained in the image of the original determinant map (since $p_1$ is general). The lemma follows.
     \end{proof}
     Let us end this section with a couple of remarks.
     \begin{rem}
    
         Using \autoref{lemma:gaussmap} instead of \autoref{rank} in the proof of \autoref{estimate}, one can obtain better bounds on the genus provided that the normal sheaf $\mathcal P$ fits in the hypothesis of \autoref{lemma:gaussmap}. This strategy will be used in the next section for curves lying on abelian varieties. 
     \end{rem}
     \begin{rem}
        The bound of \autoref{lemma:gaussmap} is sharp. For instance let us consider the following example. Pick $C\subset \mathbb
     P^{g-1}$ a trigonal canonical curve. Project the curve from a point lying on a trisecant line, the image of $C\to \mathbb P^{g-2}$ has a triple point at the image of the trisecant line. Now it is easy to see that the $\emph{kernel}$  $\mathcal P^\vee$ of this map has $g$ co-sections.
     Indeed we have an inclusion $\mathbb P^{g-2}\subset \mathbb PH^0(\mathcal P)$, now the section $s$ corresponding to the triple point vanishes in the intersection of the trisecant line with the curve in $\mathbb P^{g-1}$, hence it induces an inclusion $A\to \mathcal P$, where $A=\mathcal O_C(P_1+P_2+P_3)$ is a $g^1_3$ by geometric Riemann-Roch. Inside $\mathbb P H^0(\mathcal P)$ the linear subspaces $\mathbb P H^0(A)$ and $\mathbb P^{g-2}$ meets in one point along the class of the section $[s]$. This gives $h^0(\mathcal P)\ge g$, \autoref{lemma:gaussmap} gives $h^0(\mathcal P)\le g.$
     \end{rem}
 \section{An application to families of Jacobians with a large isotrivial factor}
 \label{section: curves isoztivial}

\quad In this section everything is assumed to be over $\mathbb C$. As an application of \autoref{estimate} we will study families of curves with a large isotrivial factor in the family of Jacobians. Using the fact that the \emph{kernel} bundle of the Gauss map is the normal bundle of the curves in the abelian variety we will use \autoref{lemma:gaussmap} to obtain better bounds than \autoref{estimate}. Recall that for a curve $D\subset A$ the Gauss map $\mc{G}$ can be simply defined as follows: let 
\begin{eqnarray*}
   \mathcal{G}^{\circ}:\tilde D^{\circ} &\to& \mathbb PT_0A \\
   \mathcal G^{\circ}(p) &= & [T_p D],
\end{eqnarray*}
where $\tilde D$ denotes a desingularization of $D$, and $D^{\circ}\subset D$ denotes the open dense locus of points $p\in {D}$ where the map $D\rightarrow A$ is an immersion. We write $\mc{G}: \tilde{D}\rightarrow \mathbb PT_0A$ for the unique map extending $\mc{G}^{\circ}$. The Gauss map is non degenerate if and only if the curve $D$ generates the abelian variety. Before going into the proof let us recall some notation from the proof of \autoref{main} and a useful preliminary lemma to bound the degree of some particular  \emph{kernel} bundles. 
\begin{setup}
 Suppose that $B_0\subset \mathrm{Hilb}_A$ is a reduced variety parametrizing irreducible reduced curves is an abelian variety $A$. Let $f: \mc D\to B_0$ be a resolution of the universal family over $B_0$ with  general member $D_b$ and let $\pi:\mathcal D \to A$ be the tautological map. Locally at $b\in B_0$ we get a commutative diagram

\[
    \begin{tikzcd}
    0 \arrow{r} & \mathcal T_{ D_b}\arrow{r}& {\mathcal T_{\mathcal D}}_{|_{D_b}}\arrow{r} \arrow{d}{d\pi}& \mathcal N_{ D_b/\mathcal D}=T_bB_0\otimes \mathcal O_{D_b}\arrow{d}{d\pi}\arrow{r} & 0 \\
    &  & {\pi^*\mathcal T_{A}}_{|_{ D_b}}\arrow{r}& \mathcal P \arrow{r} & 0,
    \end{tikzcd}
    \]
     where  $\mathcal P $ is the torsion free part of the cokernel of $d\pi:\mathcal T_{D_b}\to{\pi^*\mathcal T_{A}}_{|_{\pi(D_b)}}$. Hence $\mathcal P$ is a vector bundle. Note that, by definition, $\mathcal P^\vee$ is the \emph{kernel bundle}\footnote{We say that a vector bundle is the \emph{kernel bundle} of a morphism (or a rational map) $X\to \mathbb PV$ if it is the kernel of the evaluation $V^\vee\otimes \mathcal O_X\to \mathcal O_X(1)$.} of the Gauss map $\mathcal{G}: D\rightarrow \mathbb{P}T_0A$. i.e. its dual  fits into the exact sequence
     \[
     0\rightarrow \mathcal{P}^{\vee}\rightarrow (T_0A)^{\vee}\otimes \mathcal{O}_D\rightarrow \omega_D. 
     \]

     Indeed $T_0A^\vee=H^0(\Omega^1_A) \subset H^0(\omega_D)$ is the linear system inducing the Gauss map. This follows from the fact that given an hyperplane $H\subset \mathbb PT_0A$ there is a unique (up to scalar) one form $\eta$ such that $H$ is the projectivization of $\mathrm{ker}(\eta)$, now it is easy to check that $Z(\eta_{|_D})$ is the same as $H\cap D$ outside of the base locus of the Gauss map.


\end{setup}

\begin{lem}\label{lemma:degree3}
    Let $f:\mathcal{D}\to B_0$ be a  family of smooth curves of genus $g$ which is maximally non-isotrivial (i.e. the induced map $B_0\rightarrow \mc{M}_g$ is finite onto its image), with isotrivial factor $A_f$ in its Jacobian of dimension $q_f\geq 3$.  Suppose the following conditions hold:
    \begin{enumerate}
        \item for $b\in B_0$, the induced map $D_b\rightarrow A_f$ is birational onto its image,
        \item for all $p\in D_b$, $s\in H^0(\mc{P}\otimes \mathcal I_p)$, $\deg Z(s)\geq 3$. 
    \end{enumerate}
Then $2g-2\geq 3(\dim B_0+1)$.
    
\end{lem}

\begin{proof}
Note that $T_bB$ is disjoint from $T_0A$ inside $H^0(\mc{P})$ for $b\in B$ general, by our assumption. Hence $T_bB\oplus T_0A\subset H^0(\mathcal P)$. 

For a pair $(U, \mc{E})$ where 
\begin{itemize}
    \item $\mc{E}$ is a vector bundle on $D_b$, 
    \item $U\subset H^0(\mc{E})$,
\end{itemize}
we say that it is \emph{great} if for all $s\in U$, either $Z(s)$ is empty, or $\deg Z(s)\geq 3$. Note that the pair $(H^0(\mathcal P), \mc{P})$ is great.

    Take a section $s_1\in H^0(\mc{P})$ which does not vanish on $D_b$. Consider the quotients $\mc{P}/\mc{O}_{D_b}\cdot s_1$ and $H^0(\mathcal P)/\langle s_1 \rangle$. It follows that $(H^0(\mathcal P)/\langle s_1 \rangle, \mc{P}/\mc{O}_{D_b}\cdot s_1)$ is also great. 

    We now iterate this to produce sections 
    \[
    s_i\in H^0(\mathcal P)/\langle s_1, \cdots , s_{i-1} \rangle \subset  H^0(\mc{P}/\mc{O}_{D_b}\langle s_1, \cdots s_{i-1} \rangle)
    \]
    such that the pair 
    \[
    (H^0(\mathcal P)/\langle s_1, \cdots , s_{i} \rangle , \mc{P}/\mc{O}_{D_b}\langle s_1, \cdots s_{i} \rangle )
    \]
    is great, up to $i=q_f-2$. When $i=q_f-2$, we obtain a great pair $(W, L)$ where $L$ is a line bundle. 

    \begin{claim}\label{claim:degree3}
        The map $\varphi_W: D_b \rightarrow \mb{P}W^{\vee}$ has degree at least three onto its image.
    \end{claim}

We first show how to finish the proof of the lemma, assuming Claim \ref{claim:degree3}: indeed, 
\begin{align}
  \deg \omega_{D_b}&\geq \deg(L) \\
  &\geq  3 (\dim W-1)\\
  & \ge  3(\dim B_0-(q_f-2)+q_f-1)
\end{align}
from which we deduce the desired inequality.

\begin{proof}[Proof of Claim \ref{claim:degree3}]

For convenience, we first introduce a piece of notation. For a pair $(U, \mc{E})$ where $\mc{E}$ is a vector bundle on $D_b$ and $U\subset H^0(\mc{E})$ a sub vector space, and $\mathscr{D}$ a $\mathbb{Z}_{\geq 0}$-divisor on $D_b$, write $U(-\mathscr{D})\subset U$ for the subspace of sections vanishing along $\mathscr{D}$. 

 It now suffices to prove that  $\mathrm{Bl}(W(-p))$ is of degree $\ge 3$ for any $p$, which in turn is implied by the claim that  $\mathrm{Bl}|H^0(\mathcal{P}\otimes \mathcal I_p)|=\bigcap_{t\in H^0(\mathcal P\otimes \mathcal I_p)} Z(t)$ is of degree $\ge 3$ for any $p$. Indeed if we have a section $s\in W(-p) \subset H^0(L\otimes \mathcal I_p)$ which is not vanishing at a point $q\in \mathrm{Bl}|H^0(\mathcal P\otimes \mathcal I_p)|-\{p\}$, then  we can lift it to $\tilde s \in  H^0(\mathcal P)$ with $\tilde s$ vanishing at $p$, but of course $\tilde s$ cannot vanish at $q$, which is a contradiction because $q \in \mathrm{Bl}|H^0(\mathcal P\otimes \mathcal I_p)|$. Therefore $W(-p)=W(-\mathrm{Bl}|H^0(\mathcal P\otimes \mathcal I_p)|)$. This implies that $\varphi_W$ must contract all such base loci (i.e. $\varphi_W(\mathrm{Bl}|H^0(\mathcal P\otimes \mathcal I_p)|$ is a point for any $p$), hence it must be of degree at least $3$.

 So we are left to prove that $\mathrm{Bl}|H^0(\mathcal P\otimes \mathcal I_p)|$ is of degree $\ge 3$ for any $p$.

 Suppose $s, t\in H^0(\mathcal P\otimes \mc I_p)$, with $s$ the unique section of $T_0A_f(-p)$. Consider a complement $U$ of $\langle s \rangle$ in $T_0A_f$. Moreover, 
 \begin{itemize}
     \item let $s_q\in U$ be the unique non-zero section, up to scaling, which vanishes at $q$, 

 \item 
 and for $\lambda \in \mb{C}$ general, let $s_q^{\lambda} \in U\oplus \langle s+\lambda t\rangle$ be the unique non-zero section (up to scaling) which vanishes at $q$.
 \end{itemize}

Now assume by way of contradiction that $Z(s+\lambda t)$ has fixed part of degree $\le 2$ (hence it has moving part). This implies that $s\wedge t=0$.  Hence we get a commutative diagram (for general $q$)
\[
\begin{tikzcd}
  \langle s,t \rangle \otimes \mathcal O_{D_b}
  \arrow[drr, bend left]
  \arrow[ddr, bend right]
   & & \\
    & \mc{O}_{D_b} \arrow[r, "s_q"]  \arrow[d, "s_q^\lambda"] \arrow{dr}{s_q^U}
      & \mc{P} \arrow[d] \\
& \mc{P} \arrow[r] &\mc{P}_U, 
\end{tikzcd}
\]
 which corresponds to the geometric picture
 \[
 \begin{tikzcd}
     \mathbb PT_0A \arrow[dr,dashed,"\pi_{s}"] & D_b\arrow{r}{\mathcal{G}_{\lambda}} \arrow{d} \arrow{l}{\mathcal G} & \mathbb P( U \oplus \langle s+\lambda t \rangle ) \arrow[dl,dashed,"\pi_{\lambda}"] \\
     & \mathbb P U &.
 \end{tikzcd}
 \]
 Here $\mathcal G, \mathcal G_{\lambda}$ are the Gauss map and a perturbation of the Gauss map which both have kernel bundle $\mathcal P$; whereas $C\to \mathbb P U$ is a projection of the Gauss map with kernel bundle $\mathcal P_U$\footnote{The perturbation of the Gauss map being induced by the linear system $\bigwedge ^{q_f-1}( U \oplus \langle s+\lambda t \rangle ) \simeq ( U \oplus \langle s+\lambda t \rangle ) ^\vee \subset H^0( \omega_{D_b})$, in other words we have an exact sequence $0\to \mathcal P^\vee \to ( U \oplus \langle s+\lambda t \rangle )^\vee \otimes \mathcal O_{D_b}\to \omega_{D_b} $. Note that this is well defined only for general $\lambda$.}, fitting into a short exact sequence \[
 \langle s,t \rangle \otimes \mathcal O_C \to \mathcal P \to \mathcal P_U \to 0.
 \]
 The section $s_q^U$ is the unique section (up to scaling) of $U(-q)\subset H^0(\mathcal P_U(-q))$. Now for general $q$ by the first diagram we get $Z(s_q^{\lambda})\subset Z(s_q^U)$ (which is a finite set for general $q$)\footnote{Note that this is  defined only for general $q$, because for some $q$ the map $s_q^U$ can be zero}. This implies in particular that $Z_q=\bigcap_\lambda Z(s_q^\lambda)$ is of degree $\ge 3$ for general $q$ (they all vanish in at least $3$ points and the union of their zero loci is finite).
 Now taking a sequence  $q_i \to p$ we get that $\mathrm{lim}_{i}Z(s_{q_i}^\lambda)\subset Z(s_p^\lambda)=Z(s+\lambda t)$.
 Here the left hand side contains $\mathrm{lim}_iZ_{q_i}$ which is of degree at least $3$ and does not depend on $\lambda$, contradicting our assumption that $Z(s+\lambda t)$ has fixed part of degree $\le 2$.
\end{proof}

\end{proof}

\begin{rem}
    The bound of \autoref{lemma:degree3} can be stated in more general terms, saying that if $\mathcal P$ is a \emph{kernel} bundle of rank $\ge 2$ with the property that any section $s\in H^0(\mathcal P)$ has zero locus of degree $0$ or of degree $\ge d$ then $\mathrm{deg}(\mathcal P)\ge d(h^0(\mathcal P)-\mathrm{rk}(\mathcal P))$. This bound is sharp, for instance if $A$ is a base point free line bundle of degree $d$ with $h^0(A)=2$ then $\mathcal P=A^{\oplus r}$ fits in this hypothesis and $\mathrm{deg}(\mathcal P)=d(h^0(\mathcal P)-\mathrm{rk}(\mathcal P))$. 
    \label{degreeez}
\end{rem}
This yields the following for curves lying on abelian varieties:
\begin{thm}
    \label{genuscurvesinabelian}
    Let $A$ be an abelian variety of dimension $q_f$ and consider $B_0\subset \mathrm{Hilb}_{A_f}$ parametrizing irreducible reduced curves of geometric genus $g$ with $B_0\to \mathcal M_g$ generically finite onto its image, then
\[g\geq 
\begin{cases}
 \mathrm{dim}(B_0)+q_f  &\textrm{if}    \quad\mathrm{dim}(B_0)\le \frac{1}{2}(q_f-1)\quad \mathrm{or} \quad  2q_f-5\le \dim B_0 \\
 \frac{3 \mathrm{max}\{\dim B_0+1, q_f-1\}}{2} +1   &\textrm{if} \quad \, \,\frac{1}{2}(q_f-1) \le  \mathrm{dim}(B_0) \le  2q_f-5 
\end{cases}
\]
\label{curvesinab}
\end{thm}
\begin{proof}Let $D_b$ be the normalization of a general curve in the family.
    Consider the kernel bundle $\mathcal P$ of the Gauss map, and the exact sequence
    \[\begin{tikzcd}
        0 \arrow{r} & \mathcal P^\vee  \arrow{r} & T_0A \otimes \mathcal O_{D_b}\arrow{r} & \omega_{D_b}.
    \end{tikzcd}
    \]

 We get an inclusion $T_bB_0\oplus T_0A\subset H^0(\mathcal P)$.

 Hence an induced contraction map \[
\bigwedge^{q_f-1} H^0(\mathcal P)\to H^0(\mathrm{det}\mathcal P)\subset H^0(\omega_{D_b}).
\]
Since hyperelliptic curves are rigid in abelian varieties (see \cite[Theorem 1]{pirola_kummer}) we may suppose $D_b$ non-hyperelliptic.
Let us consider two cases
\begin{itemize}
    \item If we can find a section $s\in H^0(\mathcal P)$ vanishing at one or two points we get that $Z(s)$ is rigid hence if we consider $W\subset T_0A$ of codimension $1$ and the sequence \[
0 \to \mathcal P \to (W\oplus \langle s\rangle )^\vee \otimes \mathcal{O}_{D_b} \to \omega_{D_b},
\]
we get that $\mathcal P$ is as in the hypothesis of \autoref{lemma:gaussmap} hence we get $g\ge h^0(\mathcal P)\ge \mathrm{dim}B_0+q_f$.

\item If all the sections of $\mathcal P$ are vanishing in either $0$ points or $\geq 3$ points we get, in view of \autoref{lemma:degree3}, $2g-2\ge 3(\dim B_0+1)$. Notice we also have the inequality $2g-2\ge 3(q_f-1)$ since if we consider general sections $s_1,\dots, s_{q_f-1}$ vanishing in $q_f-1$ general points they must vanish in $3(q_f-1)$ general points by our assumptions.
\end{itemize}
So putting everything together 
\[
g\geq \min \left\{\dim(B_0)+q_f, \frac{3 \mathrm{max}\{\dim B_0+1, q_f-1\}}{2} +1 \right\},
\]
and now a tedious computation finishes the proof of the theorem. For example, if $2q_f-5\leq \dim(B_0)$, then 
\[
\frac{3 \mathrm{max}\{\dim B_0+1, q_f-1\}}{2} +1\geq \dim(B_0)+q_f,
\]
and hence $g\geq \dim(B_0)+q_f$, which proves the theorem in this case. The cases $\dim(B_0)\leq \frac{1}{2}(q_f-1)$ and $\frac{1}{2}(q_f-1) \le  \mathrm{dim}(B_0) \le  2q_f-5 $ are similar.
\end{proof}

Now we can deduce
\begin{cor}
    Let $f:\mathcal{C}\to B$ be a  family of curves of genus $g$ over $\mathbb C$, whose map to $\mc{M}_g$ is finite onto its image. Let $A_f=\mathrm{ker}( \mathrm{Alb}(\mathcal{C})\to \mathrm{Alb}{(B)})$ of dimension $q_f\ge 2$ be the isotrivial isogeny factor of the associated family of Jacobians. Then,
\[g\geq 
\begin{cases}
 \mathrm{dim}(B)+q_f  &\textrm{if}    \quad\mathrm{dim}(B)\le \frac{1}{2}(q_f-1)\quad \mathrm{or} \quad  2q_f-5\le \dim B \le 2q_f-2 \\
 \frac{3 \mathrm{max}\{\dim B+1, q_f-1\}}{2} +1   &\textrm{if} \quad \, \,\frac{1}{2}(q_f-1) \le  \mathrm{dim}(B) \le  2q_f-5 \\
 2q_f-1+\frac{\dim B}{2}&\textrm{if}  \quad \, \, 2q_f-2\leq \mathrm{dim}(B).
\end{cases}
\]
Moreover, suppose 
$\mathrm{dim}(B)>2q_f-2$, then the above bound is sharp and the only family realizing this bound is the family of degree $2$ covers of a curve of genus $q_f$ branched over $n=\mathrm{dim}(B)$ points (hence $B$ is a generically finite cover of $C^{(n)}$ the $n$-fold symmetric product of $C$ and $A_f$ is isogeneous to $\mathrm{Jac}(C)$).
\end{cor}

\begin{proof}
Notice the the image of a general curve in the family generates $A=A_f$.
We get a map $B\dashrightarrow \mathrm{Hilb}_{A}$ (we consider the image of a general curve in the family inside $A$ with the reduced structure). Let $B_0$ be the image of this  map and let $F$ be the general fiber, so that we have $\mathrm{dim}(B)=\mathrm{dim}(B_0)+\mathrm{dim}(F)$. Let us notice that in case $\mathrm{dim}(F)=0$ the theorem follows from \autoref{curvesinab}. So from now on we shall assume that the map $C_b\to A_f$ is not birational onto its image- Let $D_b=\mathrm{Im}(C_b)\subset A_f$.

By \autoref{curvesinab} the geometric genus $g_{B_0}$ of a general curve in $B_0$ (i.e. the image of a curve of the original family in the abelian variety) satisfies
\[
g_{B_0}\ge \begin{cases}
 \mathrm{dim}(B_0)+q_f  &\textrm{if}    \quad\mathrm{dim}(B_0)\le \frac{1}{2}(q_f-1)\quad \mathrm{or} \quad  2q_f-5\le \dim B_0 \\
 \frac{3 \mathrm{max}\{\dim B_0+1, q_f-1\}}{2} +1   &\textrm{if} \quad \, \,\frac{1}{2}(q_f-1) \le  \mathrm{dim}(B_0) \le  2q_f-5 
\end{cases}
\]
Now the theorem follows from using Hurwitz formula to bound $g$ in terms of $g_{B_0}$ and the dimension of $F$. Since the curves over $F$ have constant image in the abelian varieties the dimension of $F$ is bounded by the degree of the ramification divisor $C_b \to D_b\subset A_f$. By Hurwitz formula
\[
2g-2\ge 2(2g_{B_0}-2)+\mathrm{deg}(R)\ge  2(2g_{B_0}-2)+\mathrm{dim}(F).
\]
This yields $g\ge 2g_{B_0}-1+\frac{\mathrm{dim}(F)}{2}$. The  theorem now follows by  case by case computation: we illustrate this in one case and leave the rest to the reader:
\begin{itemize}
    \item If $\mathrm{dim}(B)\le \frac{1}{2}(q_f-1)$ we deduce $g_{B_0}\ge
 \mathrm{dim}(B_0)+q_f$. Then we get  \begin{eqnarray*}     
g&\ge& 2g_{B_0}-1+\frac{\mathrm{dim}(F)}{2}\\ &\ge& 2\mathrm{dim}(B_0)+2q_f-1+\frac{\mathrm{dim}(B)-\mathrm{dim}{(B_0)}}{2}\\ &\ge& \frac{3}{2}\mathrm{dim}(B_0)+\frac{3}{2}\mathrm{dim}(B)+q_f\\
&\ge& \mathrm{dim}(B)+q_f,
 \end{eqnarray*}
 as required.
\end{itemize}

\end{proof}

\section{An application to curves on general hypersurfaces}
\quad In this section everything is assumed to be over $\mathbb C$. Let us recall briefly some notation from the introduction \begin{eqnarray*}
    M_{g,Y}^{bir}&:={\{[C] \in \mathcal M_g \ | \ \exists f:C \to Y \text{ such that $f$ is birational onto its image}\}}\subset \mathcal M_g.\end{eqnarray*}
The goal of this section is to prove the following. 
\begin{thm}\label{curves vs hypersf}
     Let $Y\subset \mathbb P^{n+1}_\mathbb{C}$ be a very general hypersurface of degree $\ge 2n-1$ (if $n\ge 4$) or of degree $\ge 2n+1$ (if $n=2,3$), then  \[\dim M_{g,Y}^{bir} \le g-(d-2n).\]
\end{thm}
\quad The cases $d=2n,2n-1$ of \autoref{hypint} can be deduced in a similar way using \cite[Theorem 0.4]{Voisinhypersf} to rule out the existence of rational curves in very general hypersurfaces of degree $2n$, and to prove that any surface inside a very general hypersurface of degree $2n-1$ carries a canonical form.

\quad With a similar proof we will also deduce.
\begin{thm}\label{hypersf}
Let  $Y \subset \mathbb P^{n+1}$ be a very general hypersurface of degree $e\ge 2n+2$. Then, for a very general hypersurface $X \subset \mathbb P^{m+1}$, any rational map $X \dashrightarrow Y$ is constant.
\end{thm}

\quad Before proving \autoref{curves vs hypersf} let us compare it with the main result of Clemens in \cite{clemens}.
\begin{thm}[Clemens]
Let $\pi:\mathcal C\to B$ be a family of immersed curves on a generic hypersurface $Y\subset \mathbb P^{n+1}$ of degree $d$ and let $Z\subset Y$ be the subvariety they span. Then\[
\mathrm{codim}(Z)\ge \frac{2-2g}{\mathrm{deg}(C_b\to \mathbb P^{n+1})}+d-(n+2).
\]
\label{cle:curves vs hyp}    
\end{thm}
\quad The first thing to remark is the assumption that the curves $C_b\to Y$ are immersed (i.e. that $C_b$ is smooth and the differential of $C_b\to Y$ is nowhere vanishing), this poses constraints on the singularities of the curves. For instance the image of $C_b$ cannot have any cusp. In order to deduce bounds on the dimension of $M_{g,Y}^{bir}$ from \autoref{cle:curves vs hyp} one needs to prove that the general elements of a family of curves of genus $g$ on a generic hypersurface is an immersed curve and to be able to have a lower bound on $\mathrm{deg}(C_b\to Y)$ for curves moving in a large family in $\mathcal M_g$. Roughly speaking, this second point amounts to control the dimension of arbitrary Brill-Noether loci in $\mathcal{M}_g$ which is a difficult task. There are many components not having the expected dimension, see for instance the introduction of \cite{pflu}; notable examples are the Castelnuovo curves \cite{cili}. For example the locus of hyperelliptic curves is always a component of the locus of curves admitting a map of degree $\ge {2n+2}$ in $\mathbb P^{n+1}$ whose dimension grows asymptotically as $2g$, while the dimension of the locus of Castelnuovo curves in $\mathbb P^{n+1}$ grows asymptotically as $g+\sqrt{g}$ (c.f. \cite{cili}). This just to point out that even if one is able to carry  these two steps one may end up with a worst estimate than \autoref{curves vs hypersf} (of course the locus of hyperelliptic curves itself will not cause any trouble since it cannot lie in $M_{g,Y}^{bir}$). 

\quad The proof of \autoref{curves vs hypersf} is a consequence of \autoref{main} and results of Ein and Voisin on the positivity of the canonical bundle of subvarieties of very general hypersurfaces.
\begin{setup}\label{set up hypersurfaces} We start with some preliminaries: we will construct a subvariety $B$ of some Hilbert scheme of $Y$, up to base change, which parametrizes curves of geometric genus $g$ and having dimension $\ge \dim M_{g,Y}^{bir}$. 
Fix a general $Y$ as in \autoref{curves vs hypersf} and $g\ge 2$. Consider an irreducible component $M\subset M_{g,Y}^{bir}$ of maximal dimension. By the properties of Hilbert schemes there exist countably many quasi-projective varieties $\mathcal B_k$ parametrizing pairs $([C],f)$ where $[C]\in M$ is a genus $g$ curve and $f :C \to Y$ is birational onto its image. In fact such a map can be identified with its graph $\Gamma_f \subset C \times Y$ which is a smooth curve moving in $\mathcal C\times Y$, where $\mathcal C \to \mathcal M_g$ is the universal family. 
Consider the universal family $\mathcal C \to \mathcal M_g$ and let $\mathcal C_M \to M$ be its restriction to $M$. For each $k$, let $\pi_k: \mathcal X_k \dashrightarrow \mathcal B_k$ be the family obtained pulling-back $\mathcal C_M \to M$ via the projection $\mathcal B_k \to M$. By definition of $\mathcal B_k$ as a Hilbert scheme, there is a universal rational map
\[
\Phi_k: \mathcal X_k \dashrightarrow Y.
\]
By construction $\mathcal B_{k_0}\to M$ must be dominant, for some $k_0$ (i.e. there exists $k_0$ such that $\overline{\mathrm{Im}(\mathcal B_{k_0})}$ is not a proper closed subset). 
After slicing if needed $\mathcal B_{k_0}$ by a suitable (general) linear section and restricting our attention to an irreducible component, we get to a family $p: \mathcal X \to  \mathcal B$ of curves of genus $g$ over an irreducible base, which pulls back from $\mathcal C_M \to M$ via a generically finite dominant morphism $\mathcal B \to M$ and having a rational map $\Phi: \mathcal X \dashrightarrow Y$ which is non constant on the general fiber $X_b$ of $p$. 

\quad Let now $T\subset Y$ be the Zariski closure of $\Phi(\mathcal X)$. 
Take a desingularization $T'\to T$ and let $\mathcal X':= \mathcal X \times_T T'$ be the base change. Note there is a natural map onto some Hilbert scheme $\mu: \mathcal B \dashrightarrow \text{Hilb}_{T'}$. Such a map must be generically finite onto its image. 
So we arrived to the commutative diagram
\[
\xymatrix{
\mathcal X' \ar@{-->}[r] \ar[d] & \mathcal H_B \ar[d] \\
\mathcal B \ar@{-->}[r] &  B
}
\]
where $B\subset \text{Hilb}_{T'}$ to be the closure of $\text{im} \ \mu$ and $\mathcal H_B$ is the pull-back to $B$ of the universal family over $\text{Hilb}_{T'}$. Let $\pi: \mathcal H_B \to Y$ be the natural map. Finally, set $t:= \dim T$.
\end{setup}
\begin{rem}
    For the case $g=1$, replace in the construction $\mathcal M_1$ with $\mathcal M_{1,k}$ for some $k>0$, so that we have a universal family $\mathcal C \to \mathcal M_{1,k}$. In the end $\mu: \mathcal B \dashrightarrow B \subset \text{Hilb}_{T'}$ will have fibers of dimension $k$, which is $\dim B=1$.
\end{rem}
\quad Before proving \autoref{curves vs hypersf} and \autoref{hypersf} let us recall the definition of birationally $k-$very ample, see \cite[Definition 1.1]{BDELU} (this notion was already studied in \cite{knu1},\cite{knu2} under a different name). 
\begin{defn}
A line bundle $L$ on a projective variety $X$ \textit{satisfies the property} $(\text{BVA})_k$ if there is a Zariski-closed set, $Z(L)\subsetneq X$ depending on $L$, such that the map
\[
H^0(X,L) \to H^0(X, L \otimes \mathcal O_\xi)
\]
is a surjection for any $\xi \subset X$ subscheme of length $k+1$ whose support is disjoint from $Z(L)$.
\end{defn}
\quad We are now ready to prove \autoref{curves vs hypersf}.

\begin{proof}
    Let $\mathcal H \to B$ be the universal family. Applying \autoref{main} we get, for $b \in B$
    \[
    g \ge h^0({\omega_{T'}}_{|_{\pi (H_b)}})+ \dim B -(t-1)
    \]
    from which if $g>1$, recalling that $\dim B= \dim  M^{bir}_{g,Y}$, we have  
    \[
    g \ge  h^0({\omega_{T'}}_{|_{\pi (H_b)}}) + \dim M_{g,Y}^{bir} -(t-1).
    \]
 Results of Ein \cite{Einhypersf} and Voisin \cite{Voisinhypersf} imply that $\omega_{T'}$ satisfies $(\text{BVA})_{d+t-2n-2}$: for a proof we refer the reader to \cite[Proposition 3.8]{BDELU}.
  Therefore if we choose $\pi(C_b)$ avoiding generically $Z(\omega_{T'})$ we get 
     $ h^0({\omega_{T'}}_{|_{\pi (H_b)}}) \ge d+t-2n-2+1$ (for general $b$ the image $\pi(C_b)$ will avoid generically $Z(\omega_{T'})$). Plugging this in the previous inequality and recalling that $d\ge 2n+1$ we get the desired statement. If $g=1$ we have $\mathrm{dim}(B)=1>3g-3$ and we reach a contradiction.
\end{proof}
We are now ready  to prove \autoref{hypersf}. 
\begin{proof}[Proof of \autoref{hypersf}] 
    Assume the theorem fails i.e.: fix a general $Y$ as in the statement of the theorem and assume there is a Zariski-dense set of hypersurfaces $X$ of degree $d$ in $\mathbb P^{m+1}$ admitting a non-constant map $X \dashrightarrow Y$. We want to reduce to the case $m=1$. In order to do that, notice that the general degree $d$ plane curve is obtaining by slicing with planes a general hypersurface $X \subset \mathbb P^{m+1}$ of degree $d$; moreover if $\phi:X \dashrightarrow Y$ is non-constant, the restriction of $\phi$ to a general plane section of $X$ is still non-constant. This gives a non constant mapping from a general plane curve of degree $d$ to $Y$.
    Turning to details, let $\mathbb G$ be the Grassmanian of projective planes in $\mathbb P^{m+1}$, fix $P_0 \in \mathbb G$, $z_1, \dots z_4 \in P_0$ in general position and $z_5, \dots z_{m+3} \in \mathbb P^{m+1}$ in general position too. Let also
    \[
    \mathcal G:=\{(P,x_1, \dots x_4) \in \mathbb G \times ( \mathbb P^{m+1})^4 \ | \ x_1, \dots, x_4 \in P\}
    \]
    be the set of $4$-pointed planes.
    As in \ref{set up hypersurfaces} there are countably many schemes $\mathcal S_k, \mathcal B_k$, where $\mathcal S_k$'s parametize pairs $([X],\phi)$, with $X \in \mathbb PH^0(\mathbb P^{m+1}, \mathcal O_{\mathbb P^{m+1}}(d))$ and $\phi: X \dashrightarrow Y$ non-constant and the $\mathcal B_k$'s parametrize pairs $([C],\psi)$, with $C \in \mathbb P H^0(P_0, \mathcal O_{P_0}(d))$ and $\psi : C \to Y$ is non-constant. (Note that the $\mathcal B_k$'s could be empty a priori).
    There is a map
    \[
     \begin{aligned}
    &\pi: \bigcup_k \mathcal S_k\times \mathcal G \dashrightarrow \bigcup_k \mathcal B_k,\\
    &(X, \phi, P, x_1, \dots, x_4) \mapsto (T^{-1}(X.P), \phi \circ T),
    \end{aligned}
    \]
     where $T$ is the only transformation in $\text{PGL}(m+2)$ mapping $z_i$ to $x_i$ for $i \le 4$ (and hence mapping $P_0$ to $P$) and fixing $z_i$ for $i>4$. Moreover $\text{pr}_1 \circ \pi$ is dominant, since the general degree $d$ plane curve is a plane section of some hypersurface $X \in \mathbb PH^0(\mathbb P^{m+1}, \mathcal O_{\mathbb P^{m+1}}(d))$.
    This shows that if the theorem fails for some $m \ge 1$, it fails also for $m=1$. We thus continue the proof assuming $m=1$. We can also assume $d\ge 4$, since the case $d=3$ follows by \ref{curves vs hypersf}. The same construction done in \ref{set up hypersurfaces} applies to this case. Therefore, there exist desingularization $T'$ of some $t$-dimensional $T\subset Y$ with a subvariety of some Hilbert scheme $B \subset \text{Hilb}_{T'}$ such that the universal family $\mathcal H_B \to B$ dominates $T'$ and $B$ dominates the moduli space
    \[
   \mathcal M= \mathcal M(d,1):= \mathbb P H^0(\mathbb P^{2}, \mathcal O_{\mathbb P^{2}}(d))/\text{PGL}(3).
    \]
    \autoref{main} then yields
      \[
    \frac{(d-1)(d-2)}2 \ge h^0(\omega_{T'})-h^0(\omega_{T'} \otimes \mathcal I_{\pi (H_{b})})+ \frac{(d+2)(d+1)}2-9 -(t-1).
    \]
    As in the proof of \autoref{curves vs hypersf} we have $h^0(\omega_{T'})-h^0(\omega_{T'} \otimes \mathcal I_{\pi (H_{b})}) \ge e+t-2n-2+1$. Plugging this in the previous inequality we get $d \le 2$ (if $d=3$ one has to be slightly more careful since the above moduli space does not parametrize isomorphism classes of elliptic curves, but elliptic curves together with a degree $3$ line bundle, hence the RHS should be replaced with RHS-1, which does not affect the conclusion). Since all maps from a rational curve to $Y$ are constant (again by \cite{Einhypersf}, \cite{Voisinhypersf} for instance), this concludes  the proof.
\end{proof}
\begin{rem}
    With the same technique of \ref{curves vs hypersf} we can get information on the analogue loci for hypersurfaces. In particular, let us denote for $m,d \ge 0$,
    \[
    \mathcal M(d,m):=\mathbb PH^0(\mathbb P^{m+1}, \mathcal O_{\mathbb P^{m+1}}(d))/\text{PGL}(m+2)
    \]
    the moduli space of degree $d$ hypersurfaces in $\mathbb P^{m+1}$. Let $d,m,e,n \ge 0$ with $d \ge 2m+2, e \ge 2n+2$, fix a general $[Y]\in \mathcal M (e,n)$ and consider the locus
    \[
    M_Y^{bir}(d,m):=\{[X] \in \mathcal M(d,m) \ | \ \exists f: X \to Y \text{ such that $f$ is birational onto its image}\} \subset \mathcal M(d,m).
    \]
    Then as in \ref{curves vs hypersf} we get
    \[
    \dim M_Y^{bir}(d,m)\le {d-1 \choose m+1}-(e-2n).
    \]
\end{rem}

\bibliography{refer}
\bibliographystyle{alphaspecial}

\end{document}